\documentclass{article}
\usepackage{graphicx} 
\usepackage{amsmath}
\usepackage{amsthm}
\usepackage{hyperref}
 \usepackage{booktabs} 
 \usepackage{makecell}
 \usepackage{float}

\usepackage{subcaption}

\usepackage{amsfonts}
\usepackage{url}

\usepackage{array} 
\usepackage{makecell} 
\usepackage{cite}
\newtheorem{lemma}{Lemma}
\newtheorem{theorem}{Theorem}

\title{Homotopies and Maps between Eigenvalues of some Generalized Lucas Sequences and the Mandelbrot set}
\author{Arturo Ortiz-Tapia}
\date{May 2025}

\begin{document}

\maketitle
\begin{abstract}
The eigenvalues of companion matrices associated with generalized Lucas sequences, denoted as \( \mathcal{L} \), exhibit a striking geometric resemblance to the Mandelbrot set \( \mathcal{M} \). This work investigates this connection by analyzing the statistical distribution of eigenvalues and constructing a variety of homotopies that map different regions of \( \mathcal{L} \) to structurally corresponding subsets of \( \mathcal{M} \).

In particular, we explore both global and piecewise homotopies, including a sinusoidal interpolation targeting the main cardioid and localized deformations aligned with the periodic bulbs. We also study a variation of the Jungreis map to better capture angular and radial structures. In addition to visual and geometric matching, we classify the eigenvalues according to their dynamical behavior—identifying subsets associated with hyperbolic, parabolic, Misurewicz, and Siegel disk points.

Our findings suggest that meaningful correspondences between \( \mathcal{L} \) and \( \mathcal{M} \) must integrate both geometric deformation and dynamical classification. In light of these observations, we also suggest a conjectural homeomorphism between \( \mathcal{L} \) and a dense subset of the Mandelbrot cardioid boundary, based on the behavior of the sinusoidal homotopy and the eigenvalue accumulation.
\par
Finally, we prove that the sinusoidal homotopy defines a homeomorphism (modulo countable exceptions) from \( \mathcal{L} \) onto the boundary of the Mandelbrot cardioid, and, when composed with Douady's tuning map, extends to any baby cardioid or stable region of \( \mathcal{M} \), reinforcing the structural correspondence between these sets.
\end{abstract}

\section{Introduction}

The Mandelbrot set ($\mathcal{M}$ hereinafter) is an important starting point for exploring dynamical systems, and many systematic relations between points in the complex plane (or the complex Riemann sphere) and $\mathcal{M}$ have been deemed worthy of exploration \cite{devaney1999mandelbrot}. 

\par
In a previous work, it was found that the union of sets of eigenvalues of matrices corresponding to certain generalized Lucas sequences, denoted hereinafter by $\mathcal{L}$, produced a graphical representation resembling a slightly enlarged version of $\mathcal{M}$, though not an exact scaled copy. This work aims to conduct a statistical analysis of these eigenvalues, providing insights into suitable numerical explorations of maps and homotopy-like transformations between $\mathcal{L}$ and $\mathcal{M}$, of the following sort: 
If $C \subseteq \mathbb{R}^n$ is a \emph{convex} subset of Euclidean space and $f, g: [0, 1] \to C$ are paths with the same endpoints, then there is a \textbf{linear homotopy}\cite{Hatcher2002} (or \textbf{straight-line homotopy}) given by 
\begin{align}
  H: [0, 1] \times [0, 1] &\longrightarrow C \\\nonumber
                   (s, t) &\longmapsto (1 - t)f(s) + tg(s).  \\\nonumber  
\end{align}

\par
 It is important to notice that $\mathcal{L}$ is countably infinite, and that thus, homotopically transformed into a correspondingly countably infinite subset of points belonging to, or related with (by geometrical resemblace) $\mathcal{M}$, the details of which will be explained in each section.

\par
The structure of this paper is as follows. The first part expands on the findings of the previous work, focusing on statistical properties of the eigenvalues and developing intuition about potential homotopies applicable to their structure. The next section introduces the Jungreis map, which was initially defined for the complement of the unit circle but extends naturally to include it. Applying this map to the eigenvalue set $\mathcal{L}$ reveals that its boundary maps onto $\mathcal{M}$, while the remaining points correspond to regions escaping to infinity at different "speeds." 

\par
Following this, we explore a transformation that takes all points to the \emph{cardioid} (the period-1 component of $\mathcal{M}$). The final key section introduces a conjectured \emph{generalized piecewise homotopy}, where selected subsets of $\mathcal{L}$ are mapped to specific regions of $\mathcal{M}$. Notably, the \emph{antennas of $\mathcal{L}$} are mapped to the corresponding \emph{bulbs} of $\mathcal{M}$, such as those of \emph{period 6 over period 4}, while most other points transform into the \emph{cardioid}. Points that are part of the \emph{period-2 bulb} (behind the cardioid) are mapped accordingly, except for the \emph{tail of $\mathcal{L}$}, which is mapped to the \emph{tail of $\mathcal{M}$}. Since the spacing of eigenvalues in the tail of $\mathcal{L}$ remains to be analyzed, determining the final transformation to $\mathcal{M}$ requires further refinement. This piecewise homotopy thus consists of several maps working together to deform $\mathcal{L}$ into $\mathcal{M}$.
\par
As part of this construction, we also provide a formal proof that this homotopy constitutes a homeomorphism (modulo countable exceptions) from \( \mathcal{L} \) to the boundary of the Mandelbrot cardioid, and—via Douady's tuning map—to the boundaries of baby Mandelbrot sets and other stable regions of \( \mathcal{M} \), anchoring our visual and numerical observations in topological terms.

\par
The final section presents a global discussion of these findings, synthesizing the results and addressing broader implications. Throughout this work, the term "map" may be used interchangeably to refer to functions, homotopies, or other types of transformations, unless further clarification is required for readability.

\par

\section{Statistics of \texorpdfstring{$|\lambda|$}{|lambda|}}\label{Sec:Intervals-Statistics}

\begin{figure}[ht]
\centering
\includegraphics[width=5in]{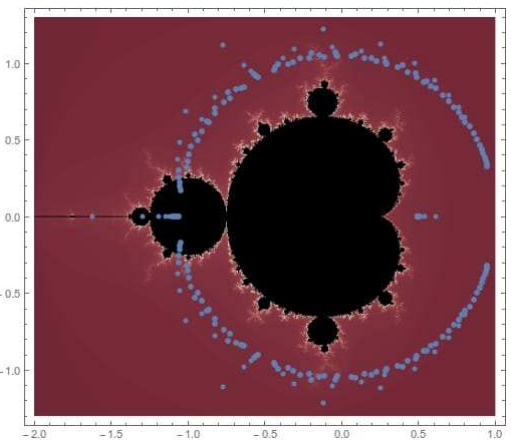} 
\caption{A plot of $\mathcal{M}$ and the reciprocal of the eigenvalues coming from $A_{200}$.}
\label{fig:MandelbrotMyEigenvalues}
\end{figure}

In a previous work \cite{tapia2018golden}, the union of sets of arithmetical reciprocal of the eigenvalues of these matrices ($1/\lambda_n$ with $n$ being the size of the companion matrix), so, $1/\lambda_n\in \bigcup_{k=1}^n\det (A_n-I\lambda)=\mathcal{L}$) rendered a plot that geometrically resembled $\mathcal{M}$ (Fig.\ref{fig:MandelbrotMyEigenvalues}). This compelling graphical resemblance motivated this work, to attempt to use several types of maps, in order to actually formalize a transformation that takes those $1/\lambda_{n}$ to a point on the boundary of $\mathcal{M}$ (Fig. \ref{fig:MandelbrotMyEigenvalues}).

\pagebreak

\subsection{Matrix Construction}
Let us brief recapitulate some aspects of the previous work, which will be used in this and other sections.

We construct a series of matrices $A_n$ of increasing size $n$, where

\begin{equation}
n \in \{10, 20, 50, 100, 150, 200, 250, 300, 350, 400, 450, 500\}.    
\end{equation}

Each \emph{companion} matrix $A_n$ has the following structure \cite{tapia2018golden}:

\begin{equation}\label{Eq:CompanionMatrix}
A_n = \begin{pmatrix}
1 & 1 & 1 & \cdots & 1 \\
1 & 0 & 0 & \cdots & 0 \\
0 & 1 & 0 & \cdots & 0 \\
\vdots & \vdots & \ddots & \ddots & \vdots \\
0 & 0 & \cdots & 1 & 0
\end{pmatrix}
\end{equation}

or symbolically:
\begin{eqnarray}\label{Eq:SymbolicCompanionMatrices}
a_{1,j} &=& 1, \, 1 \leq j \leq n, \\\nonumber
a_{2,1} &=& 1, \\\nonumber
a_{j+2,j+1} &=& 1, \, 1 \leq j \leq n-2. \\\nonumber
\end{eqnarray}

\subsection{The unit circle}
For each matrix size $n$, we calculate the percentage of eigenvalues that fall within this $\epsilon$-neighborhood of the unit circle. Let $\lambda_n$ be the set of eigenvalues for matrix $A_n$, and $\lambda_n^\epsilon$ be the subset of eigenvalues satisfying the proximity criterion. 

The results of our numerical experiment (Table \ref{tab:unit_circle_distribution}) provide insights into the distribution of eigenvalues for the companion matrices. We observe that as the matrix size $n$ increases, the trend shows that most eigenvalues fall within the the $\epsilon$-neighborhood of the unit circle. However, we notice that there are still many eigenvalues left out.

\begin{table}[ht]
\centering
\begin{tabular}{|c|p{2.5cm}|p{2.5cm}|p{2.5cm}|p{2.5cm}|}
\hline
\textbf{n} & \% \textbf{\makecell{within \\ Unit Circle}} & \textbf{\makecell{\# Eigen \\ Solutions}} & \textbf{\makecell{\# Within \\ Unit Circle}} & \textbf{\makecell{\# Left \\ Out}} \\ \hline
10 & 0.0 & 54 & 0 & 54 \\ \hline
20 & 4.78469 & 209 & 10 & 199 \\ \hline
50 & 18.3673 & 1274 & 234 & 1040 \\ \hline
100 & 40.6813 & 5049 & 2054 & 2995 \\ \hline
150 & 71.3882 & 11324 & 8084 & 3240 \\ \hline
200 & 83.631 & 20099 & 16809 & 3290 \\ \hline
250 & 89.3542 & 31374 & 28034 & 3340 \\ \hline
300 & 92.4915 & 45149 & 41759 & 3390 \\ \hline
350 & 94.3996 & 61424 & 57984 & 3440 \\ \hline
400 & 95.6483 & 80199 & 76709 & 3490 \\ \hline
450 & 96.5114 & 101474 & 97934 & 3540 \\ \hline
500 & 97.1337 & 125249 & 121659 & 3590 \\ \hline
\end{tabular}
\caption{Behavior of eigenvalue distribution: percentage within the unit circle and the count of solutions left out as \(n\) (the matrix size) increases.}
\label{tab:unit_circle_distribution}
\end{table}

We will approximate the behavior of the number of points within the unit circle with a function that has a known asymptotic behavior:\\

Let the logistic function 
\begin{equation}
L / (1 + \exp[-k \cdot (x - x_0)])
\end{equation}

with $\text{Model fit Parameters=}\{L\to 94.4725,k\to 0.0271419,\, x_0\to111.11\}$ represent the asymptotic behavior of those points, as the cumulative percentage approaches 100\%:
\begin{equation}
\lim_{n\to\infty} \frac{94.4725}{1 + e^{-0.0271419(n-111.11)}} = 94.4725 
\end{equation}

\begin{figure}[ht]
\centering
\includegraphics[width=5in]{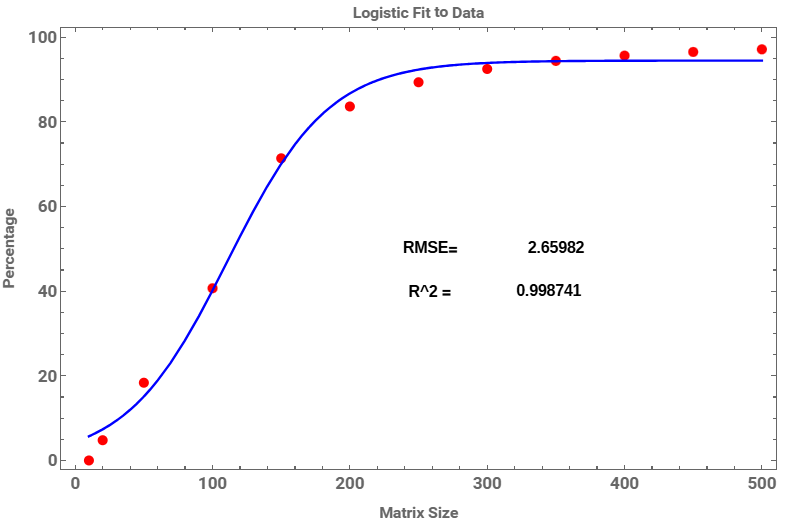} 
\caption{Fit of the logistic equation to the data, with $\text{Model fit Parameters =}\{L\to94.4725,k\to0.0271419,x_0\to 111.11\}$.}
\label{fig:figure1}
\end{figure}

As \(n \to \infty\), the percentage \( \% \text{ within Unit Circle } \) asymptotically approaches 100\%, while the count of eigenvalues left out continues to grow at each step (see Table \ref{tab:unit_circle_distribution}). This would imply that using the Jungreis map will transform some of the $|\lambda|\neq 1\pm\epsilon$ to other regions that are not at the boundary of $\mathcal{M}$, but are related nonetheless to it. We will explore this relationship in more detail in Section\ref{Sec:JungreisMap}.

To ensure that the exponential decay part (which is part of the denominator) of the logistic model passes exactly through the point where \( x = 10 \) and \( f(10) = 100\% \), we adjust the model as follows. We start with the general form of the exponential decay model:
\[
f(x) = a e^{-b x} + c,
\]
where \( a \), \( b \), and \( c \) are the parameters to be determined.Using the condition \( f(10) = 100 \), we substitute \( x = 10 \) into the equation:
\[
f(10) = a e^{-b \cdot 10} + c = 100.
\]
Rearranging for \( a \), we get:
\[
a = (100 - c) e^{10b}.
\]

Substituting \( a \) back into the original model reduces the number of free parameters to \( b \) and \( c \). The updated model becomes:
\[
f(x) = (100 - c) e^{-b (x - 10)} + c.
\]

This form ensures that the model satisfies \( f(10) = 100 \) exactly, while still allowing the parameters \( b \) and \( c \) to be fitted to the data.

Table~\ref{tab:eigenvalue_distribution} summarizes the distribution of eigenvalues in $\mathcal{L}$, grouped into interval bins for different matrix sizes, regardless of whether they are close to the unit circle or not. Recall that 

\begin{equation}
    \mathcal{L}=\bigcup_{k=2}^{n} \sigma(A_k)
\end{equation}
where $\sigma(A_k)$ denotes the spectrum (set of eigenvalues) of the matrix $A_k$

\begin{table}[htbp]
\centering
\caption{Distribution of eigenvalues of $\mathcal{L}$ in different interval bins, for various matrix sizes. Bins are of the form $[a, b)$, ensuring that each eigenvalue is uniquely classified. Some bins increase indefinitely as $n$ grows, while others stabilize, indicating structural spectral properties. Empty bins have been omitted to improve clarity.}
\label{tab:eigenvalue_distribution}
\begin{tabular}{|c|c|c|c|c|c|c|}
\hline
Interval & $n=10$ & $n=20$ & $n=100$ & $n=300$ & $n=500$ & $n=1000$ \\ \hline
[0.49, 0.51) & 6 & 16 & 96 & 296 & 496 & 996 \\ \hline
[0.51, 0.53) & 1 & 1 & 1 & 1 & 1 & 1 \\ \hline
[0.53, 0.55) & 1 & 1 & 1 & 1 & 1 & 1 \\ \hline
[0.60, 0.62) & 1 & 1 & 1 & 1 & 1 & 1 \\ \hline
[0.99, 1.01) & 0 & 10 & 2054 & 41759 & 121659 & 496409 \\ \hline
[1.01, 1.03) & 0 & 24 & 2550 & 2745 & 2745 & 2745 \\ \hline
[1.03, 1.05) & 4 & 36 & 219 & 219 & 219 & 219 \\ \hline
[1.05, 1.07) & 2 & 55 & 62 & 62 & 62 & 62 \\ \hline
[1.07, 1.09) & 6 & 25 & 25 & 25 & 25 & 25 \\ \hline
[1.09, 1.11) & 8 & 15 & 15 & 15 & 15 & 15 \\ \hline
[1.11, 1.13) & 5 & 5 & 5 & 5 & 5 & 5 \\ \hline
[1.13, 1.15) & 7 & 7 & 7 & 7 & 7 & 7 \\ \hline
[1.15, 1.17) & 4 & 4 & 4 & 4 & 4 & 4 \\ \hline
[1.19, 1.21) & 1 & 1 & 1 & 1 & 1 & 1 \\ \hline
[1.21, 1.23) & 4 & 4 & 4 & 4 & 4 & 4 \\ \hline
[1.29, 1.31) & 1 & 1 & 1 & 1 & 1 & 1 \\ \hline
[1.35, 1.37) & 2 & 2 & 2 & 2 & 2 & 2 \\ \hline
[1.60, 1.62) & 1 & 1 & 1 & 1 & 1 & 1 \\ \hline
\end{tabular}
\end{table}

Clearly, some of the bins for the intervals reach a constant and remain unchanged regardless of the increasing size of the companion matrix. Other bins, such as those for the intervals $[0.49,0.51)$ or $[0.99,1.01)$—this last essentially capturing the neighborhood of the unit circle—exhibit unbounded growth (as will be rigorously shown later). The following four graphics illustrate this growth, which has been fitted with an exponential function. 

%
\begin{figure}[htb]
\centering
\includegraphics[width=4in]{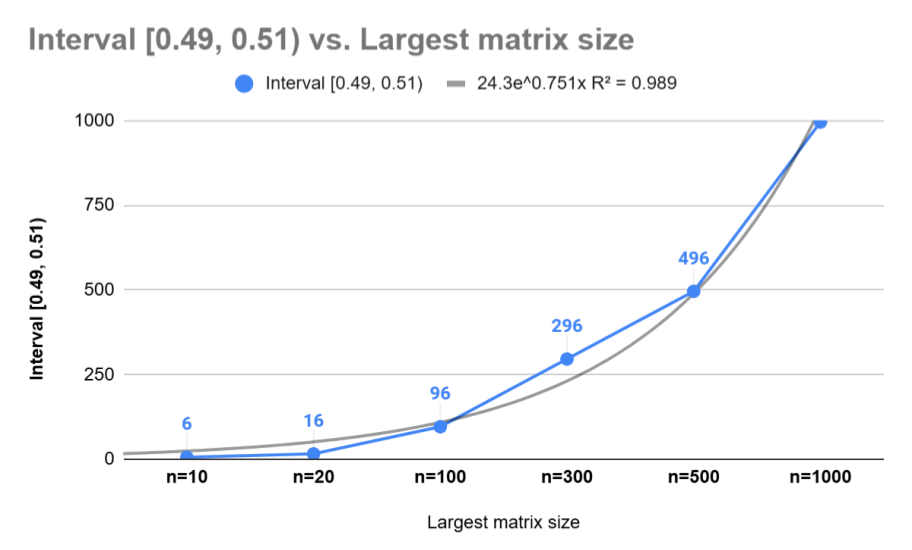} 
\caption{Bin growth for interval (radius) 0.49--0.51.}
\label{fig:Fig049_051-bin-growth}
\end{figure}
%
%
%
\begin{figure}[htb]
\centering
\includegraphics[width=4in]{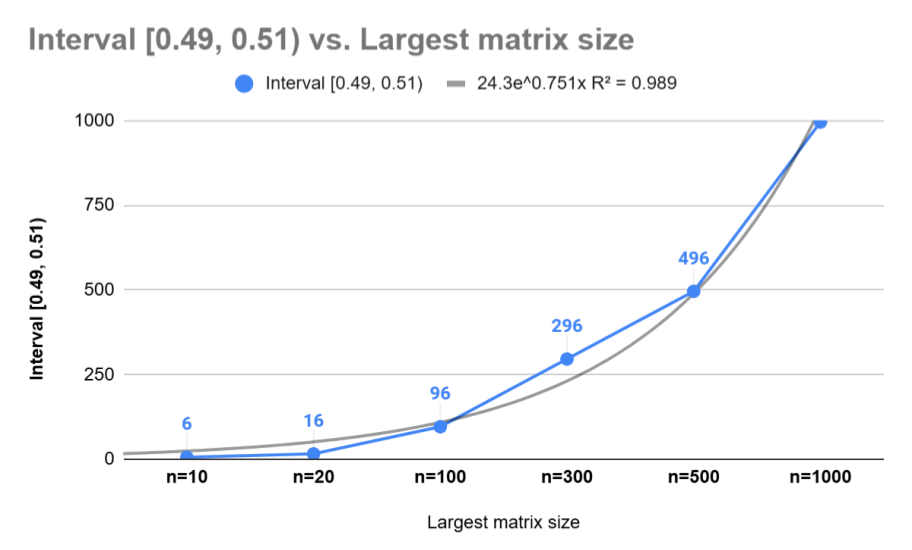} 
\caption{Bin growth for interval (radius) 0.99--1.01.}
\label{fig:Fig099-101-bin-growth}
\end{figure}
\pagebreak
%
The only bins that keep growing are those related to the boundary of the unit circle, namely $0.99-1.01$ and those below the value of the inverse of the Golden ratio, namely, the interval $[0.49, 0.51)$.
\par
Let us notice now the general pattern for the characteristic polynomial of the companion matrix:
\begin{eqnarray}\label{Eq:CharacteristicPolynomials}
    x^N-\left(\sum_{k=0}^{N-1}x^k\right),& & N\text{ is odd}\\\nonumber
    -x^N+\left(\sum_{k=0}^{N-1}x^k\right),& & N\text{ is even}\\\nonumber
\end{eqnarray}
If we set 
\begin{equation}\label{Eq:characteristicEq}
    \pm x^N\mp\left(\sum_{k=0}^{N-1}x^k\right)=0,
\end{equation}
then, given that all elements of the polynomial have a coefficient equal to $|1|$, then Eq.\ref{Eq:characteristicEq} always has at least one non-zero eigenvalue, by the fundamental theorem of algebra \cite{fine1997fundamental}.
\par

\begin{lemma} \label{lem:polyIrreducible}

The polynomials in Eq.\ref{Eq:CharacteristicPolynomials} are irreducible in the rationals.     
\end{lemma}

\begin{proof}
    Eisenstein criterion gives a \textit{sufficient} condition for a polynomial to be irreducible:If there exists a prime number \( p \) such that the following three conditions all apply\cite{eisenstein1850irreductibilitat}:

\begin{itemize}
    \item \( p \) divides each \( a_i \) for \( 0 \leq i < n \),
    \item \( p \) does \textit{not} divide \( a_n \), and
    \item \( p^2 \) does \textit{not} divide \( a_0 \),
\end{itemize}

then the polynomial is irreducible over the rational numbers.

But no such prime number number \( p \) exists, given that all coefficients are equal to $|1|$ in Eq.\ref{Eq:CharacteristicPolynomials}, and unity is neither a prime, nor any number that can be factored into prime numbers. 

The following is a \emph{necessary} condition: 
\par
\textit{A polynomial is considered irreducible over a field (like the field of rational numbers, denoted as \(\mathbb{Q}\)) if it cannot be written as the product of two non-constant polynomials with coefficients in that field.
}
\par
Since the polynomials in \textit{Eq.~\ref{Eq:CharacteristicPolynomials}} have all powers of $x$ from 0 to $N$, and since the coefficients of all those powers (including the constant term) are $|1|$, no rule can be applied to decompose those polynomials as the product of two non-constant polynomials with coefficients in \(\mathbb{Q}\)).
\end{proof}

\begin{lemma} \label{lem:density}
The set of eigenvalues of companion matrices of all sizes is dense on the unit circle in the complex plane.
\end{lemma}

\begin{proof}
For any companion matrix, its characteristic polynomial has coefficients in \( \mathbb{Q} \), ensuring that its eigenvalues are roots of a polynomial with rational coefficients. By the Fundamental Theorem of Algebra, an \( n \times n \) companion matrix has exactly \( n \) complex eigenvalues (counting multiplicity). By lemma~\ref{lem:polyIrreducible}, the characteristic polynomials of the companion matrix are irreducible; when the characteristic polynomial is irreducible over \( \mathbb{Q} \), these eigenvalues are algebraic numbers that do not belong to \( \mathbb{Q} \).

Since the set of algebraic numbers is dense in \( \mathbb{C} \), it follows that for any point \( z \) on the unit circle and any given \( \varepsilon > 0 \), we can find a sufficiently large integer \( N \) such that some \( n \times n \) companion matrix, with \( n \leq N \), has an eigenvalue \( \lambda \) satisfying \( |\lambda - z| < \varepsilon \). This establishes the density of the eigenvalues of companion matrices on the unit circle.
\end{proof}

The above lemma does not entail that all $1/\lambda$ must be in the unit circle, as it can be seen from the tables \ref{tab:eigenvalue_distribution} and \ref{tab:eigenvalue-classification}. However, it is useful to define a (piecewise) homotopy, where one of the pieces is related to the unit circle.


\section{Uniformization and the Jungreis Map}\label{Sec:JungreisMap}

To analyze the geometric structure relevant to this work, we employ the \textit{Jungreis map}, a conformal transformation that uniformizes the complement of the Mandelbrot set via a Laurent series expansion. This transformation plays a crucial role in mapping eigenvalues from the companion matrices of Lucas sequences to the dynamical structure of the Mandelbrot set.

The Jungreis map is formally given by:
\begin{equation}
    \Psi_M: \hat{\mathbb{C}} \setminus \overline{\mathbb{D}} \to \hat{\mathbb{C}} \setminus M
\end{equation}
where $\Psi_M$ is the inverse of the \textit{B\"ottcher function} $\Phi_M$, which provides a uniformizing conformal map for the complement of the Mandelbrot set. The function $\Psi_M$ can be expressed as a \textit{Laurent series}:
\begin{equation}
    c = \Psi_M (w)  =  w + \sum_{m=0}^{\infty} b_m w^{-m} = w -\frac{1}{2} + \frac{1}{8w} - \frac{1}{4w^2} + \frac{15}{128w^3} + \dots
\end{equation}
where $w \in \hat{\mathbb{C}} \setminus \overline{\mathbb{D}}$, and $c$ is a parameter in the complement of the Mandelbrot set.

Since explicit computation of the Jungreis map requires truncating this Laurent series, we focus on the \textit{numerical implementation} of this transformation. The next section details the practical method used to approximate $\Psi_M$, ensuring computational efficiency while preserving the required accuracy. Also, by "complement" of the unit disk, we have separated those values that are $|z|=1+ \epsilon$, and those that are $|z|>1+\epsilon$

\subsection{Numerical Implementation}

In order to actually use the Jungreis map, we require a formulation that can be computed numerically. Formally, let $w(t) = \cos(t) + i \sin(t)$ be the complex exponential representation of the unit circle parameterized by $t \in [0, 2\pi)$. The Jungreis map can then be approximated by truncating the Laurent series:
\begin{equation}\label{eqnJungreisMap}
J(t) = w(t) + \sum_{k=1}^{N} \frac{a_k}{w(t)^{k-1}},
\end{equation}
where $\{a_k\}$ are the coefficients of the Laurent series, determined numerically and fixed for this analysis (see Appendix~\ref{app:JungreisCoefficients} for the full set of coefficients). The resulting map $J(t)$ produces a trajectory in the complex plane as $t$ varies over its domain.

The Laurent series is truncated at $N = 65$, as this number of coefficients provides a balance between computational efficiency and sufficient approximation accuracy for the problem at hand, in particular, when applying the Jungreis map for those points with $|z| = 1 \pm \epsilon$. Since eigenvalues are calculated numerically, they are approximations, but they can still be used meaningfully for the purposes of this work.

\subsection{Behavior of Eigenvalues According to Their Modulus}

In this subsection, we analyze the behavior of eigenvalues computed from a family of $n \times n$ matrices. The matrices were constructed recursively with entries forming a structured pattern: the first row consisting entirely of ones, the first element of the second row being one, and the subdiagonal elements set to one, producing a characteristic recurrence relation. The eigenvalues were then computed numerically for matrix sizes $2 \leq n \leq 150$.

Once the eigenvalues were obtained, their multiplicative inverses were computed to facilitate a transformation analysis. The resulting eigenvalues were then expressed as coordinate pairs $(\Re \lambda, \Im \lambda)$ for visualization in the complex plane. These eigenvalues were categorized into three main sets based on their modulus:

\begin{itemize}
    \item \textbf{Eigenvalues on the unit circle}: $0.99 \leq |\lambda| \leq 1.01$
    \item \textbf{Eigenvalues inside the unit disk}: $|\lambda| < 0.99$
    \item \textbf{Eigenvalues outside the unit disk}: $|\lambda| > 1.01$
\end{itemize}

A homotopy transformation was then applied to each of these sets. The homotopy function is defined as:
\[
H(z, t) = (1 - t) z + t \left(z + \sum_{k=1}^{m} \frac{c_k}{z^{k-1}}\right),
\]
where the coefficients $c_k$ correspond to the Jungreis series expansion. This transformation allows us to observe the deformation of eigenvalues under the homotopy, particularly focusing on those near the unit circle and those with $|\lambda| > 1$.

The figures below illustrate the results of this transformation for different values of $t$. The first set corresponds to eigenvalues with $|\lambda| > 1$, while the second set corresponds to eigenvalues approximately on the unit circle.

\begin{figure}[H]
    \centering
    \begin{subfigure}[b]{0.48\textwidth}
        \includegraphics[width=\textwidth]{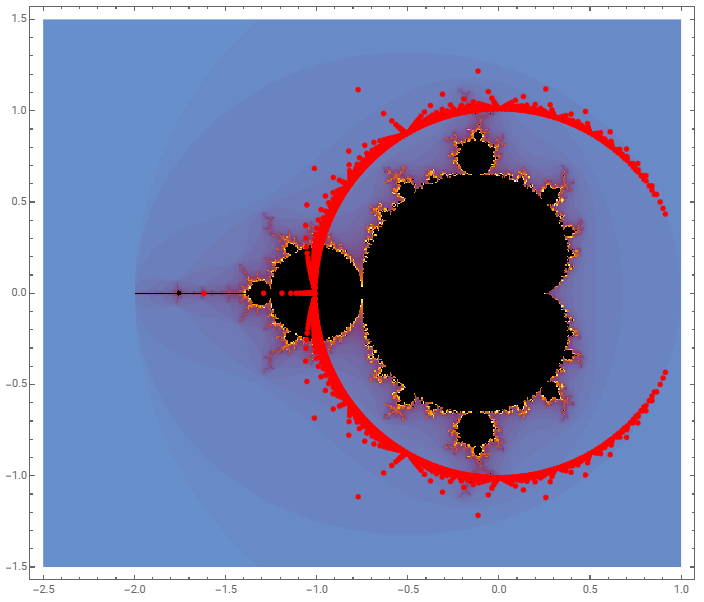}
        \caption{$t = 0$: Initial eigenvalue distribution}
    \end{subfigure}
    \begin{subfigure}[b]{0.48\textwidth}
        \includegraphics[width=\textwidth]{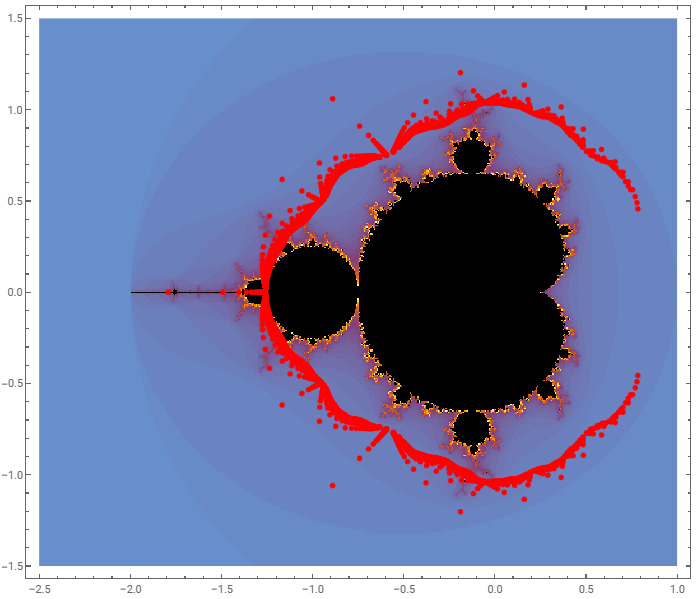}
        \caption{$t = 0.25$: Partial transformation}
    \end{subfigure}

    \begin{subfigure}[b]{0.48\textwidth}
        \includegraphics[width=\textwidth]{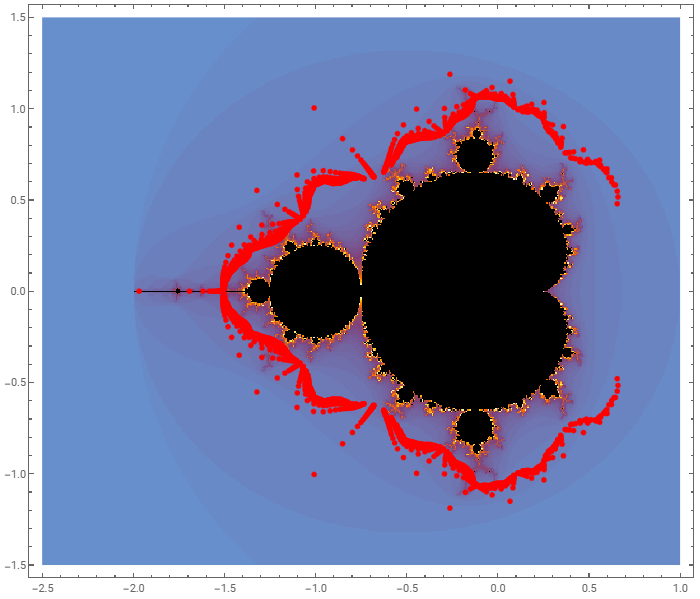}
        \caption{$t = 0.5$: Intermediate transformation}
    \end{subfigure}
    \begin{subfigure}[b]{0.48\textwidth}
        \includegraphics[width=\textwidth]{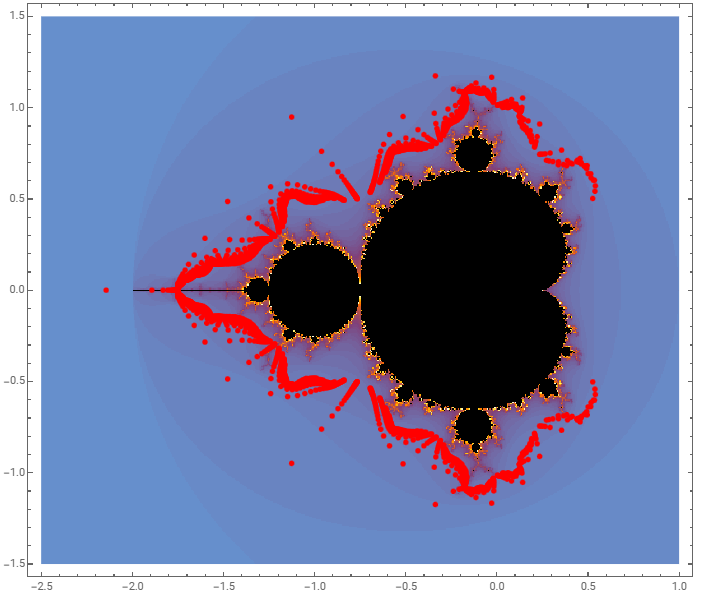}
        \caption{$t = 0.75$: Near-final transformation}
    \end{subfigure}

    \begin{subfigure}[b]{0.48\textwidth}
        \includegraphics[width=\textwidth]{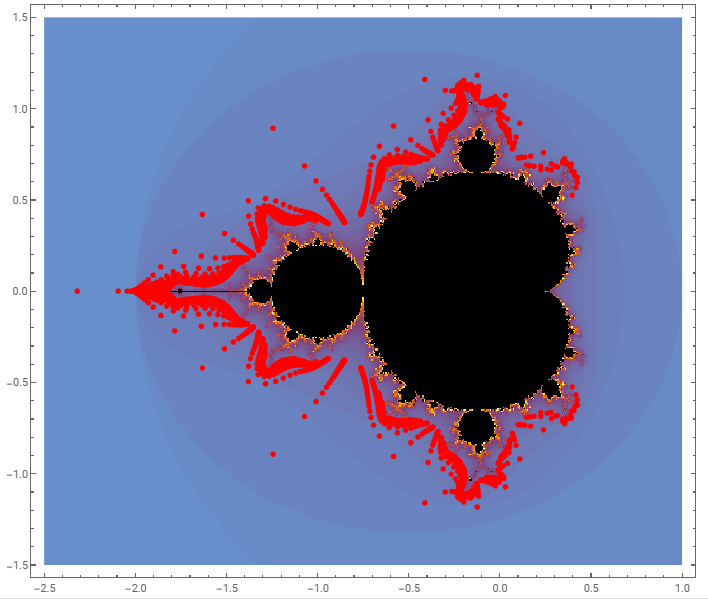}
        \caption{$t = 1$: Fully transformed distribution}
    \end{subfigure}
    \caption{Homotopy transformation for eigenvalues with $|\lambda| > 1$. The deformation is driven by the Jungreis coefficients.}
    \label{fig:homotopy_gt1}
\end{figure}

\begin{figure}[htbp]
    \centering
    \begin{subfigure}[b]{0.48\textwidth}
        \includegraphics[width=\textwidth]{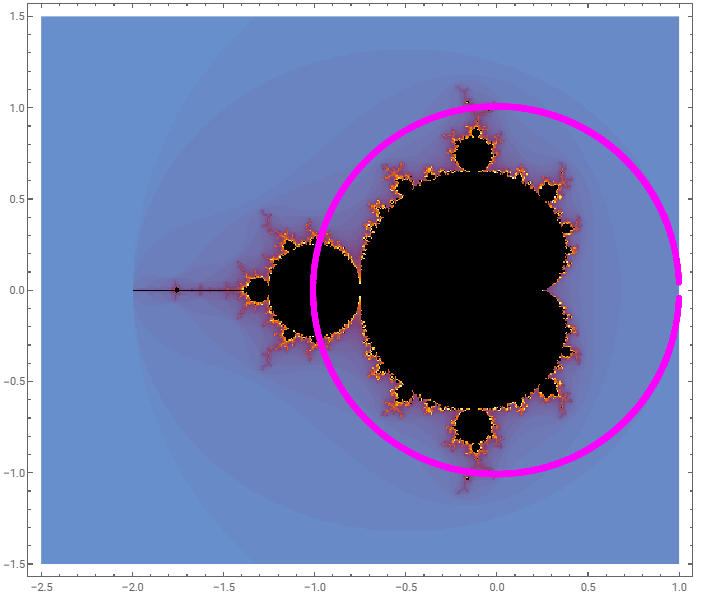}
        \caption{$t = 0$: Initial eigenvalue distribution}
    \end{subfigure}
    \begin{subfigure}[b]{0.48\textwidth}
        \includegraphics[width=\textwidth]{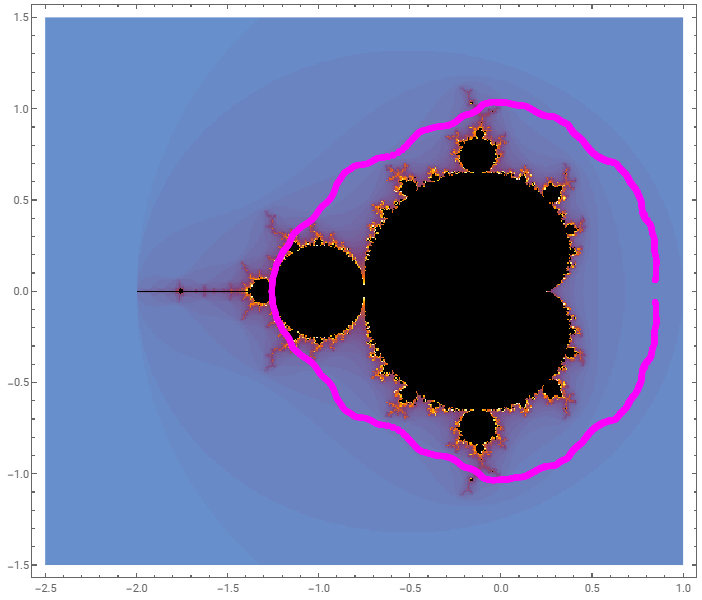}
        \caption{$t = 0.25$: Partial transformation}
    \end{subfigure}

    \begin{subfigure}[b]{0.48\textwidth}
        \includegraphics[width=\textwidth]{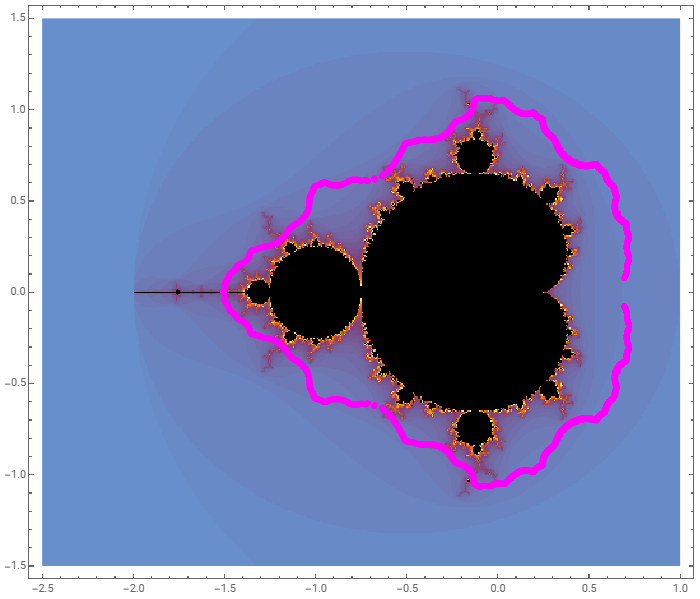}
        \caption{$t = 0.5$: Intermediate transformation}
    \end{subfigure}
    \begin{subfigure}[b]{0.48\textwidth}
        \includegraphics[width=\textwidth]{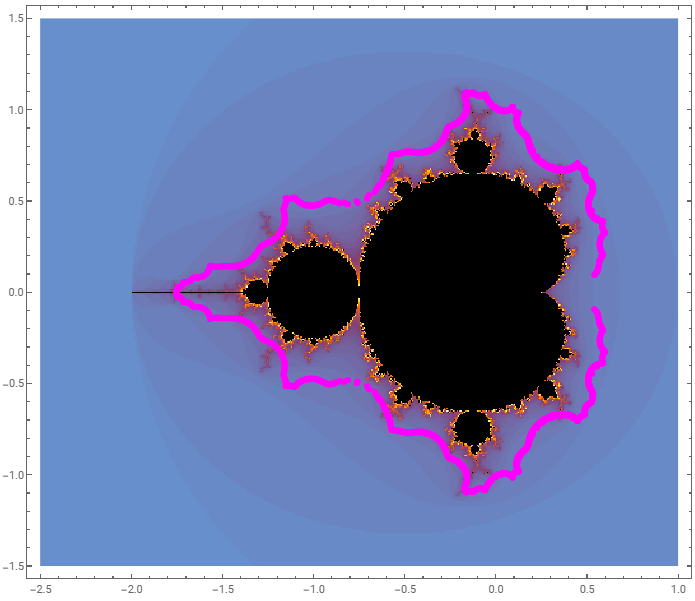}
        \caption{$t = 0.75$: Near-final transformation}
    \end{subfigure}

    \begin{subfigure}[b]{0.48\textwidth}
        \includegraphics[width=\textwidth]{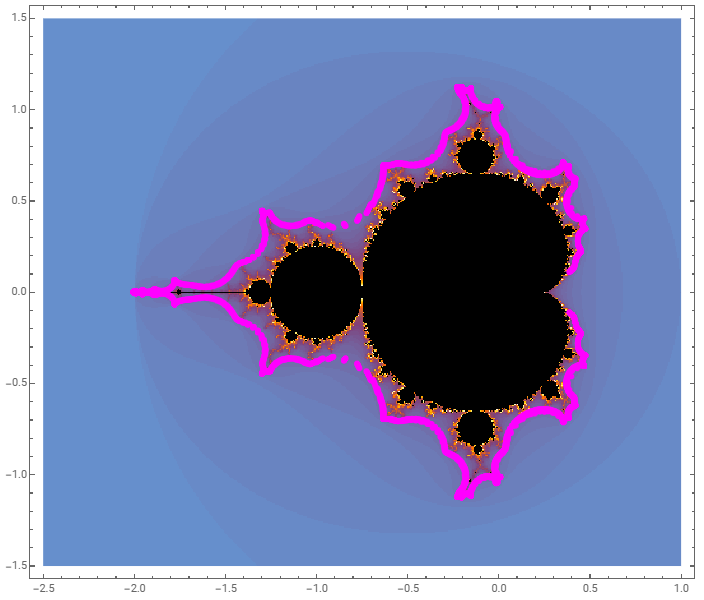}
        \caption{$t = 1$: Fully transformed distribution}
    \end{subfigure}
    \caption{Homotopy transformation for eigenvalues with $|\lambda| \approx 1$. The deformation follows the Jungreis map.}
    \label{fig:homotopy_unit}
\end{figure}

The two sets of figures illustrate the evolution of eigenvalues under the homotopy transformation, capturing the effects of Jungreis coefficients on both the outer and unit-circle eigenvalues. These transformations reveal structural properties of the eigenvalue distribution that align with the behavior predicted by the recurrence relations in the original matrix sequence.

\section{Homotopy from \texorpdfstring{$\mathcal{L}$}{L} to the cardioid of \texorpdfstring{$\mathcal{M}$}{M}}

In this section, we formalize the homotopy that transforms the eigenvalues of the companion matrices associated with the Lucas sequence towards the Mandelbrot set, specifically aligning with the boundary of the main cardioid. This transformation provides a systematic method for mapping $\mathcal{L}$ onto the cardioid of $\mathcal{M}$, reinforcing the idea that $\mathcal{L}$ can be understood as an initial structure that can be continuously deformed into a subset of $\mathcal{M}$. The homotopy is constructed to deform the eigenvalues through a piecewise function that contracts, shifts, and reshapes points in a way that gradually interpolates between their original distribution and the intended shape within the Mandelbrot set. 


Let $H: [0,1] \times \mathbb{C} \to \mathbb{C}$ be the homotopy function defined as

\begin{equation}
H(t, z) = \begin{cases}
(1-t)z + t \left( \frac{1}{4} \left( 2 \cos\theta - \cos 2\theta \right) + \frac{i}{4} \left( 2 \sin\theta - \sin 2\theta \right) \right), & |z| > 0.5, \\
(1 - 0.5t)z - t, & \text{otherwise}.
\end{cases}
\end{equation}

where $\theta = \arg(z)$ is the angular coordinate of $z$ in polar form.

The first case in the piecewise function handles the deformation of the circular part of the eigenvalue distribution, interpolating it towards the analytic boundary of the Mandelbrot set's main cardioid. The interpolation utilizes the parameterized representation of the cardioid given by \cite[p. 5]{yates1947handbook}

\begin{equation}
C(\theta) = \frac{1}{4} (2 \cos \theta - \cos 2\theta) + \frac{i}{4} (2 \sin \theta - \sin 2\theta),
\end{equation}

which can also be derived from the conformal parameterization\cite{BrownChurchill2009} of the Mandelbrot set's main cardioid using the transformation $c = \frac{z}{4} (2 - z)$\cite{Rossmannek2017}.

The second case in the homotopy function performs a general shrinking and shifting operation, moving points towards the complex plane region where the Mandelbrot set resides while preserving relative positioning. The parameter $t$ varies from 0 to 1, smoothly transitioning the eigenvalues from their initial positions to their transformed locations.

To verify the homotopy numerically, we compute the eigenvalues of the companion matrices associated with the Lucas sequence. 

\begin{figure}[htb]
\centering
\begin{subfigure}{0.45\textwidth}
    \includegraphics[width=\textwidth]{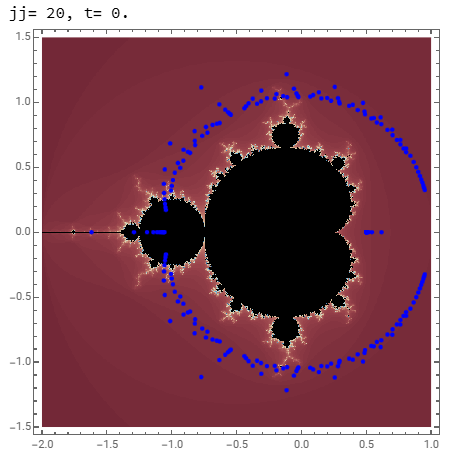}
    \caption{$H(0, z)$}\label{fig:Homotopy_t0}
\end{subfigure}
\begin{subfigure}{0.45\textwidth}
    \includegraphics[width=\textwidth]{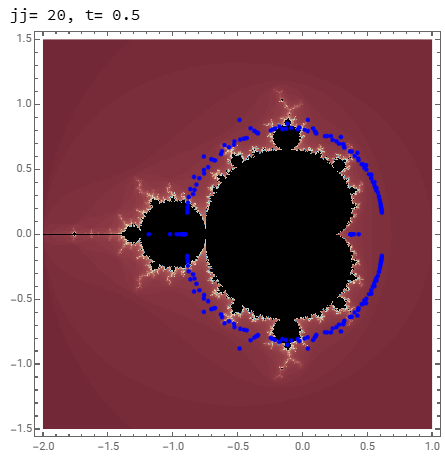}
    \caption{$H(0.5, z)$}\label{fig:Homotopy_t0_5}
\end{subfigure}

\begin{subfigure}{0.45\textwidth}
    \centering
    \includegraphics[width=\textwidth]{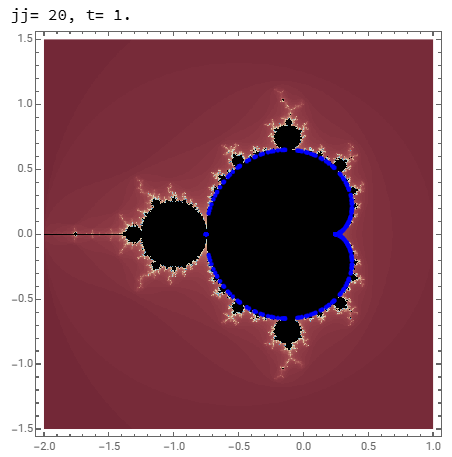}
    \caption{$H(1, z)$}\label{fig:Homotopy_t1}
\end{subfigure}
\caption{Cardioid Homotopy transition from $\mathcal{L}$ to $\mathcal{M}$.}
\label{fig:Homotopy_comparison}
\end{figure}

The numerical eigenvalues of these matrices are then mapped through $H(t, z)$ and plotted iteratively for increasing values of $t$ to visualize the deformation, as illustrated in Figs.~(\ref{fig:Homotopy_t0}, \ref{fig:Homotopy_t0_5}, and \ref{fig:Homotopy_t1}). Without loss of generality, these figures utilize only the eigenvalues from the companion matrices up to $n=20$ to optimize computational efficiency and rendering speed.

The reader may notice that all the points from the eigenset \( \mathcal{L} \) were transformed onto the cardioid of \( \mathcal{M} \), as opposed to the other homotopies shown in this work, which use different piecewise functions depending on distinguishing features within subsets of \( \mathcal{L} \). Another important remark is that the existence of self-affine, scaled-down, and rotated copies of \( \mathcal{M} \) (often called baby Mandelbrot sets) follows from the tuning theorem, first introduced by Douady and Hubbard \cite{douady1987algorithms}.

The eigenset \( \mathcal{L} \), formed from the eigenvalues of companion matrices of generalized Lucas recursions, exhibits remarkable compatibility with stable regions of the Mandelbrot set \( \mathcal{M} \). In particular, we observe that continuous homotopies can be constructed from \( \mathcal{L} \) into both the main cardioid and selected bulbs of \( \mathcal{M} \), using a combination of Douady’s tuning map and geometrically guided deformation.

Douady’s tuning map is a well-defined holomorphic transformation that associates each point \( c \) in the main cardioid of \( \mathcal{M} \) to a corresponding point in a baby Mandelbrot set. It is given by:
\begin{equation}
    \mathcal{T}(c) = c' + \lambda (c - c_0),
\end{equation}
where \( c_0 \) is the center of the main cardioid, \( c' \) is the center of a baby Mandelbrot cardioid or bulb, and \( \lambda \) is a scaling factor ensuring self-similarity. In our numerical experiments, we used the values
\begin{equation}
    c_0 = -0.75, \quad c' = -0.1575 + 1.0325i, \quad \lambda = 0.01,
\end{equation}
which yielded excellent alignment between the image of \( \mathcal{L} \) under \( \mathcal{T} \) and the observed structure of the target region in \( \mathcal{M} \).

Following this tuning, we constructed an explicit homotopy that continuously deforms the tuned eigenset into the boundary of the target cardioid. This homotopy is parameterized by \( t \in [0, 1] \), where \( t=0 \) corresponds to the tuned image and \( t=1 \) to the final transformed state. It is defined by:
\begin{equation}
    H(t, z) = (1 - t) z + t \mathcal{H}_{\text{rot}}(\theta),
\end{equation}
where \( \theta = \arg(z - c') \) and \( \mathcal{H}_{\text{rot}}(\theta) \) denotes a cardioid deformation function rotated by \( 30^\circ \) clockwise and scaled according to
\begin{equation}
    \mathcal{H}_{\text{rot}}(\theta) = s(t) e^{-i\pi/6} \left( \frac{1}{4}(2 \cos \theta - \cos 2\theta) + \frac{i}{4}(2 \sin \theta - \sin 2\theta) \right),
\end{equation}
with
\begin{equation}
    s(t) = (1 - t) \cdot 0.01 + t \cdot 0.0055.
\end{equation}

This construction results in a close match between the transformed eigenset and the front structure of the baby Mandelbrot cardioid, particularly near \( t = 1 \). While the frontal region of the deformation aligns well with the cardioid's boundary, the rear portion of the eigenset exhibits a slight protrusion beyond the back of the cardioid, suggesting that the homotopy may preserve certain geometric aspects of the eigenset's original configuration.

\begin{figure}[H]
    \centering
    \begin{subfigure}{0.45\textwidth}
        \includegraphics[width=\textwidth]{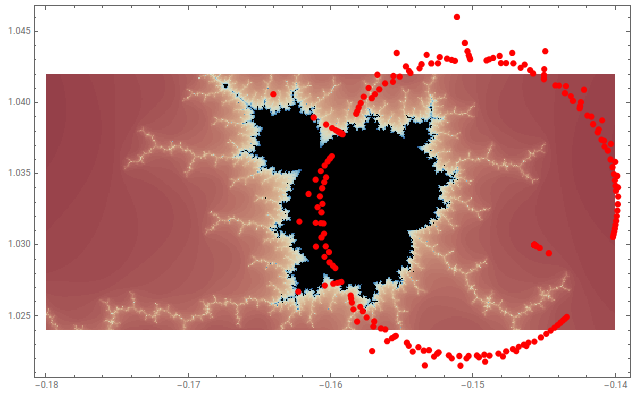}
        \caption{Homotopy at \( t=0 \): The initial transformation under tuning.}
        \label{fig:Homotopy_t0_tuning}
    \end{subfigure}
    \hfill
    \begin{subfigure}{0.45\textwidth}
        \includegraphics[width=\textwidth]{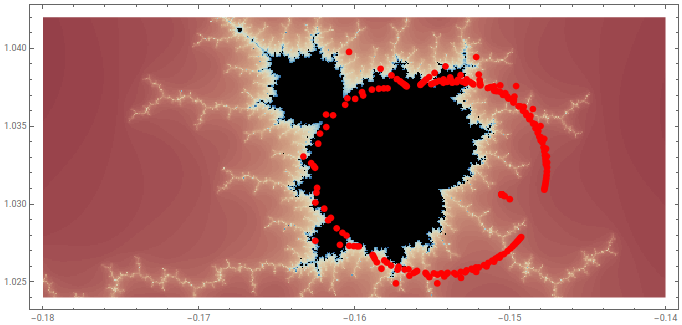}
        \caption{Homotopy at \( t=0.5 \): Intermediate transition.}
        \label{fig:Homotopy_t0_5_tuning}
    \end{subfigure}

    \begin{subfigure}{0.45\textwidth}
        \centering
        \includegraphics[width=\textwidth]{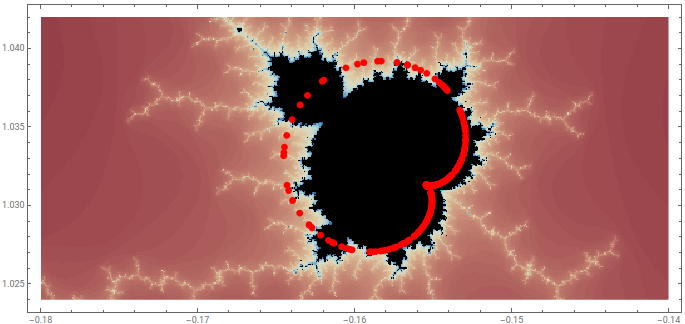}
        \caption{Homotopy at \( t=1 \): The final structure matching the baby Mandelbrot cardioid.}
        \label{fig:Homotopy_t1_tuning}
    \end{subfigure}

    \caption{Homotopy transformation from TheConstruct to the baby Mandelbrot cardioid. 
    The visualization is an approximation of a theoretically well-defined homotopy arising from the combination of Douady's tuning and a continuous deformation process.}
    \label{fig:Homotopy_tuning_cardioid}
\end{figure}

Figure~\ref{fig:Homotopy_tuning_cardioid} displays three key steps in this process. The sequence illustrates the homotopy progression, with \( \mathcal{L} \) smoothly deforming into the front structure of the baby Mandelbrot cardioid. While the frontal region of the deformation aligns well with the cardioid's boundary, the rear portion of the eigenset exhibits a slight protrusion beyond the back of the cardioid, suggesting that the homotopy may preserve aspects of the eigenset's original geometry.

\pagebreak

\vspace{1em}
\subsection{A Homeomorphism between \texorpdfstring{$\mathcal{L}$}{L} and the Boundary of the Mandelbrot Cardioid}

We now provide a formal proof that the sinusoidal homotopy introduced previously defines a homeomorphism (modulo countable exceptions) from the eigenset $\mathcal{L}$ to the boundary of the Mandelbrot cardioid $\partial\mathcal{C}$.

\begin{theorem}[Homeomorphism between $\mathcal{L}$ and the Mandelbrot Cardioid Boundary]
Let $\mathcal{L} \subset \mathbb{C}$ denote the set of multiplicative inverses of the eigenvalues of all companion matrices associated with generalized Lucas sequences. Let $\partial\mathcal{C} \subset \mathcal{M}$ be the boundary of the main cardioid of the Mandelbrot set.

Define the mapping
\[
f(z) := C(\theta) \quad \text{with} \quad \theta = \arg(z),
\]
where
\[
C(\theta) = \frac{1}{4}(2\cos \theta - \cos 2\theta) + \frac{i}{4}(2\sin \theta - \sin 2\theta),
\]
and $\arg(z)$ is taken in the principal branch $[0,2\pi)$. Partition the domain into four closed intervals:
\[
\theta \in
\begin{cases}
[0, \pi/2),\\
[\pi/2, \pi),\\
[\pi, 3\pi/2),\\
[3\pi/2, 2\pi).
\end{cases}
\]

Then the mapping $f: \mathcal{L} \to \partial\mathcal{C}$ is continuous, bijective on its image (except for countably many points), and has a continuous inverse away from those exceptional points. Thus, $f$ defines a homeomorphism from $\mathcal{L}$ onto a dense subset of $\partial\mathcal{C}$ up to a countable set of removable singularities.
\end{theorem}

\begin{proof}
\textbf{(1) Continuity via Homotopy.} Define the homotopy:
\[
H(t, z) = (1 - t) z + t \cdot C(\arg(z)), \quad t \in [0,1].
\]
This is continuous in both $t$ and $z$, provided $z \neq 0$, which holds for all $z \in \mathcal{L}$.

\textbf{(2) Injectivity on Subintervals.} On each of the four subintervals, $\theta \mapsto C(\theta)$ is injective since it locally parameterizes the cardioid boundary in a smooth and monotonic way. Restricting $\arg(z)$ accordingly prevents overlaps and ensures that $f(z) = C(\arg(z))$ is injective on $\mathcal{L}$ up to countable exceptions.

\textbf{(3) Surjectivity onto a Dense Subset.} Since $\mathcal{L}$ is a dense, countable set on the unit circle (as shown in Lemma 2), and $C(\theta)$ traces $\partial\mathcal{C}$ smoothly over $\theta \in [0, 2\pi]$, the image $f(\mathcal{L})$ is dense in $\partial\mathcal{C}$.

\textbf{(4) Inverse Continuity.} Each arc defined by a partitioned interval admits a continuous inverse map $C^{-1}$, except possibly at a countable set of boundary points. Such exceptions are removable, and the overall map remains invertible up to countably many punctures.

\textbf{(5) Extension via Douady Tuning.} Let $c_0$ be the center of the main cardioid, $c'$ the center of a finite-period bulb, and $\lambda$ a scaling factor. Then Douady’s tuning map is:
\[
\mathcal{T}(c) = c' + \lambda(c - c_0).
\]
This affine transformation is a homeomorphism. Composing $f$ with $\mathcal{T}$ yields homeomorphisms from $\mathcal{L}$ onto the boundary of any stable component or baby Mandelbrot cardioid.

\end{proof}

\begin{figure}[H]
    \centering
    \begin{subfigure}{0.45\textwidth}
        \includegraphics[width=\textwidth]{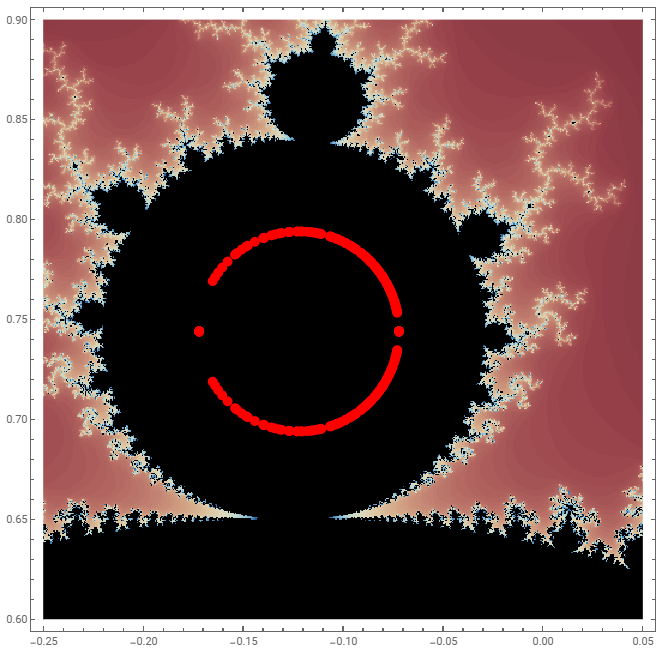}
        \caption{Homotopy at \( t=0 \): After tuning toward the period-3 region.}
        \label{fig:Homotopy_bulb3_t0}
    \end{subfigure}
    \hfill
    \begin{subfigure}{0.45\textwidth}
        \includegraphics[width=\textwidth]{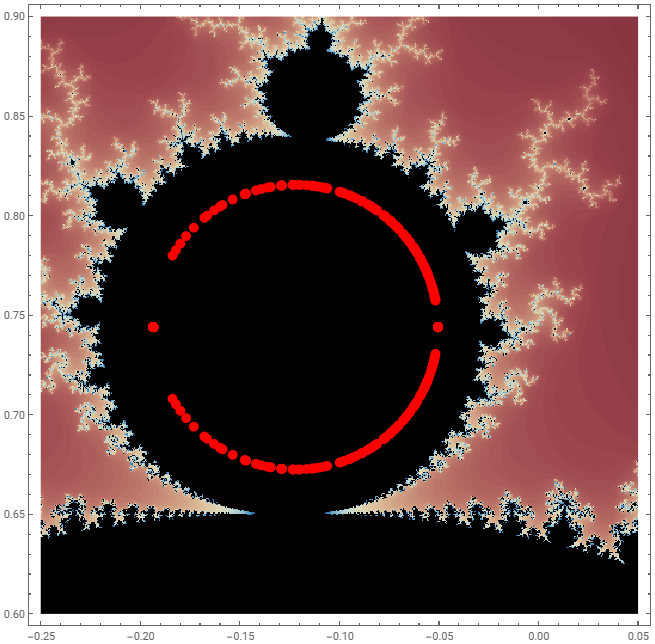}
        \caption{Homotopy at \( t=0.5 \): Expanding towards the bulb.}
        \label{fig:Homotopy_bulb3_t05}
    \end{subfigure}

    \begin{subfigure}{0.45\textwidth}
        \centering
        \includegraphics[width=\textwidth]{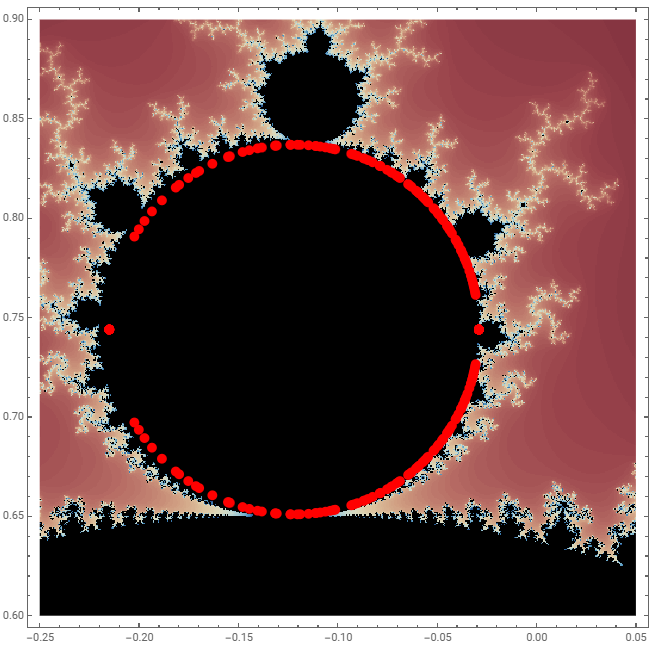}
        \caption{Homotopy at \( t=1 \): Approximate match with the period-3 bulb.}
        \label{fig:Homotopy_bulb3_t1}
    \end{subfigure}

    \caption{Homotopy from \( \mathcal{L} \) to the period-3 bulb of the Mandelbrot set. 
    The transformation expands the eigenset radially from its tuned configuration to better fit the geometry of the target bulb.}
    \label{fig:Homotopy_bulb3}
\end{figure}

\paragraph{Remark.} The adjacent figures (e.g., Fig.~\ref{fig:Homotopy_tuning_cardioid} and Fig.~\ref{fig:Homotopy_bulb3}) illustrate the intermediate and final steps of this deformation visually. These visual sequences complement the analytical proof, demonstrating geometric continuity and shape preservation during the homotopy from $\mathcal{L}$ to various cardioidal boundaries.

\par
We also explored a second mapping that leads to the top bulb of period 3. This construction again uses Douady’s tuning map to relocate the eigenset \( \mathcal{L} \) to the region surrounding a baby Mandelbrot set, this time centered at
\begin{equation}
    c' = -0.122 + 0.744i,
\end{equation}
with the zoom window
\begin{equation}
    \text{Re} \in [-0.25, 0.05], \quad \text{Im} \in [0.6, 0.9].
\end{equation}

Unlike the cardioid deformation, no polar parameterization was employed in this case. Instead, we used a radial homotopy centered at \( c' \), interpolating between initial and final radii. The deformation is defined by:
\begin{equation}
    H(t, z) = (1 - t)(r_0 e^{i\theta} + c') + t(r_1 e^{i\theta} + c'),
\end{equation}
where \( \theta = \arg(z - c') \), and in this experiment we set \( r_0 = 0.050 \) and \( r_1 = 0.093 \). 

Although this transformation does not result in a cardioid, the final image of \( \mathcal{L} \) at \( t = 1 \) shows a visually strong alignment with the top period-3 bulb. This suggests that Douady’s tuning, when combined with radial homotopies, can also be adapted to target higher-period stable components of \( \mathcal{M} \). Fig.\ref{fig:Homotopy_bulb3} illustrates this process.

\pagebreak
\section{Piecewise Homotopies}
A piecewise function is defined to address each of the parts of the eigenset; although several alternatives for defining piecewise homotopies can be taken, here we are going to explore basically two: the first is using a purely scalar transformation and the second corresponding to the cases named in a master dissertation \cite{redona1996mandelbrot}.

\subsection{Piecewise Scaling Transformation}
To construct a meaningful homotopy, we introduce a piecewise-defined transformation that adjusts scaling factors for different regions of $\mathcal{L}$. The transformation applies different weights depending on the real and imaginary parts of each eigenpoint $(x,y)$. The scaling function $S: \mathbb{R}^2 \to \mathbb{R}^2$ is defined as:
\begin{equation}
S(x,y) =
\begin{cases}
\left(\frac{0.25}{0.618} x, y\right), & x \geq 0.5, \, y = 0, \\[8pt]
\left(\frac{-1.95}{-1.6} x, y\right), & x \leq -1, \, y = 0, \\[8pt]
\left(0.7 x, \frac{0.25}{0.4} y \right), & \text{otherwise}.
\end{cases}
\end{equation}

\subsection{Homotopy Construction}
The homotopy function $H: [0,1] \times \mathbb{R}^2 \to \mathbb{R}^2$ is defined as a linear interpolation between the original points in $\mathcal{L}$ and their scaled counterparts:
\begin{equation}
H_t(x, y) = (1 - t) (x, y) + t \cdot S(x, y),
\end{equation}
where $t \in [0,1]$ represents the homotopy parameter.
\pagebreak
\subsection{Snapshots of the Homotopy}
\begin{figure}[h]
    \centering
    \begin{subfigure}[b]{0.47\textwidth}
        \includegraphics[width=\textwidth]{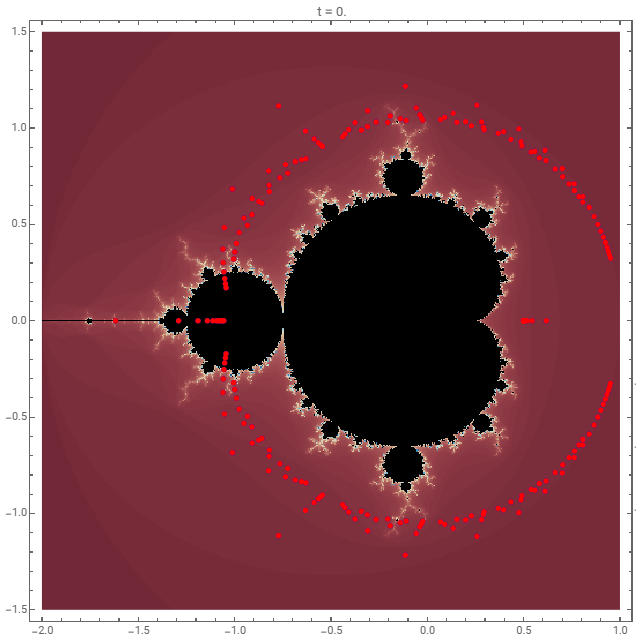}
        \caption{$t = 0$: Initial configuration $\mathcal{L}$}
        \label{fig:homotopy1}
    \end{subfigure}
    \hfill
    \begin{subfigure}[b]{0.45\textwidth}
        \includegraphics[width=\textwidth]{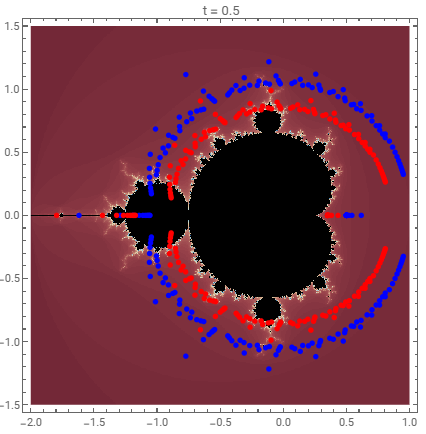}
        \caption{$t = 0.5$: Intermediate transformation}
    \end{subfigure}
    
    \vspace{0.5cm} 

    \begin{subfigure}[b]{0.45\textwidth} 
        \centering
        \includegraphics[width=\textwidth]{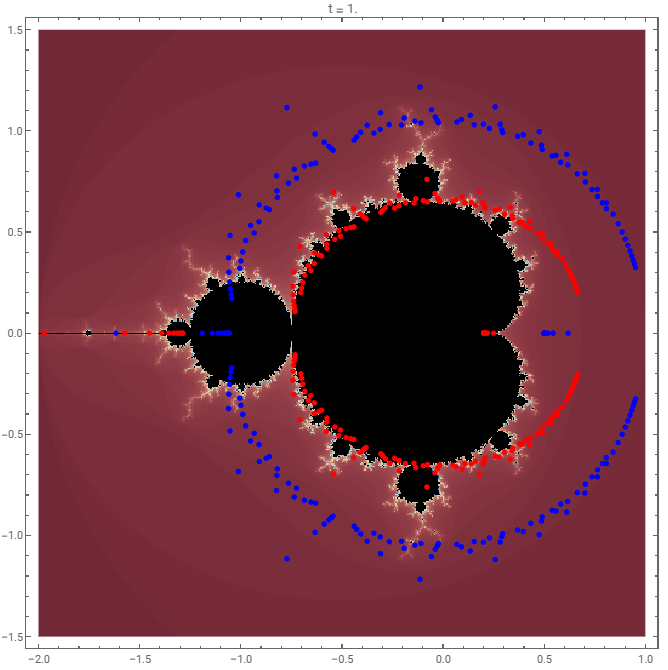}
        \caption{$t = 1$: Final transformed set}
        \label{fig:homotopy3}
    \end{subfigure}

    \caption{Evolution of the scalar, piecewise homotopy from the initial configuration $\mathcal{L}$ at $t=0$ to the final transformed set at $t=1$.}
    \label{fig:homotopy2}
\end{figure}

To visualize the transition, the homotopy is computed for discrete values of $t$:
\[
t = 0, \quad 0.25, \quad 0.5, \quad 0.75, \quad 1.
\]
Figures \ref{fig:homotopy1}, \ref{fig:homotopy2}, and \ref{fig:homotopy3} illustrate the evolution of $\mathcal{L}$ under the homotopy $H_t$.

\subsection{Observations}
It is evident from the visualization that a simple rescaling operation (i.e., multiplication by a scalar) does not necessarily lead to a meaningful transformation from $\mathcal{L}$ to $\mathcal{M}$, in part because the fundamental difference in their cardinality forces any homotopy to choose a subset of $\mathcal{M}$: $\mathcal{L}$ is a countable set, whereas $\mathcal{M}$ is uncountable. Thus, unless theoretical justifications for the selection of landing points are found, any proposed piecewise homotopy remains heuristic rather than canonical.


\subsection{Possible Homotopies that would follow from the Modulus of Eigenvalues}
In 1983, Dennis Sullivan \cite{sullivan1983conformal} showed that the following four cases describe some of the possible connected Julia sets that are associated with the M-set (the only one remaining is called a Herman ring).

\begin{itemize}
   \item \textbf{Case \#1: Hyperbolic Points.} Eigenvalues with \( |\lambda| < 1 \) are classified as hyperbolic, corresponding—by analogy—to attracting cycles in complex dynamics. In the context of the iteration \( f_c(z) = z^2 + c \), such attracting behavior occurs when the multiplier \( f_c'(z_0) \) at a fixed point \( z_0 \) satisfies \( |f_c'(z_0)| < 1 \). Parameters \( c \) for which this happens lie in the interior of the Mandelbrot set, and the associated Julia set \( J_c \) is a Jordan curve. The main cardioid of \( \mathcal{M} \) contains all parameters with an attracting fixed point, while each attached bulb corresponds to a periodic attractor of a specific period \cite{douady1984orsay}.

    \item \textbf{Case \#2: Misurewicz Points.} While the condition \( |\lambda| > 1 \) does not define Misurewicz points in the strict dynamical sense, we use it here as a proxy to identify eigenvalues whose dynamical behavior resembles repelling or preperiodic structures. In the iteration of \( f_c(z) = z^2 + c \), a Misurewicz point is a parameter \( c \) for which the critical point is strictly preperiodic, i.e., its orbit lands on a periodic cycle after finitely many steps but is not itself periodic \cite{milnor2006dynamics}. These points lie on the boundary of the Mandelbrot set and are associated with Julia sets that are dendrites—connected sets with no interior. In this classification, eigenvalues with \( |\lambda| > 1 \) are heuristically associated with such preperiodic behavior, given their tendency to "escape" under homotopies that reflect the dynamical features of \( \mathcal{M} \).
    
    \item \textbf{Case \#3: Parabolic Points.} To classify eigenvalues with \( |\lambda| \approx 1 \), we first isolate those satisfying \( 1 - \epsilon \leq |\lambda| \leq 1 + \epsilon \), for a small tolerance \( \epsilon \). Within this annulus, we heuristically classify an eigenvalue as parabolic if its argument is numerically close to zero, that is, if \( \frac{1}{\pi} \arctan\left(\frac{\text{Im}(\lambda)}{\text{Re}(\lambda)}\right) \approx 0 \). This selects eigenvalues near the positive real axis, which under homotopies may correspond to root points of hyperbolic components in the Mandelbrot set, known to exhibit parabolic dynamics. While this approach does not identify parabolic multipliers in the strict dynamical sense, it is inspired by similar angular heuristics discussed in \cite{redona1996mandelbrot}, and serves as a numerical proxy for locating eigenvalues that resemble the structural role of parabolic points in complex dynamics.

    \item \textbf{Case \#4: Siegel Disks.} The remaining eigenvalues with \( |\lambda| \approx 1 \), but whose argument does not appear close to any rational multiple of \( \pi \), are heuristically classified as Siegel-type. In complex dynamics, Siegel disks arise when a fixed point has a neutral multiplier with irrational rotation number that satisfies the Brjuno condition, resulting in quasiperiodic behavior and the formation of an invariant Jordan curve around the fixed point \cite{buff2000quadratic}. Although our classification is numerical and based solely on angular heuristics, eigenvalues with \( |\lambda| \approx 1 \) and apparently irrational arguments exhibit behavior that resembles that of Siegel disks when mapped under homotopies to the boundary of the Mandelbrot set. This interpretation follows the angular classifications used in \cite{redona1996mandelbrot} and is consistent with standard dynamical theory as presented in \cite{milnor2006dynamics}.

\end{itemize}

Beyond these four well-known cases, certain eigenvalues with \( |\lambda| \approx 1 \) do not fit neatly into either the parabolic or Siegel disk categories. Among these, \textit{Cremer points} correspond to values of \( c \) where the rotation number is irrational but fails to satisfy the Brjuno or Diophantine conditions needed for a Siegel disk. In such cases, the iteration does not lead to an invariant Jordan curve, and the Julia set remains non-locally connected with no stable periodic structure \cite{milnor2006dynamics}. 

Additionally, another exceptional case occurs when the Fatou set contains an annular region where the dynamics are conjugate to an irrational rotation on a circle. These structures, called \textit{Herman rings}, form a distinct category of Julia sets and were identified as the last remaining connected type in Sullivan's classification \cite{sullivan1983conformal}. However, Herman rings are absent in the Mandelbrot set, in contrast to Siegel disks, because the existence of a Herman ring requires a rational function of degree at least three \cite{shishikura1987quadratic}. Since the Mandelbrot set arises from iterating quadratic polynomials, Herman rings cannot occur in its associated dynamics.

So, let's start by classifying the amount of eigenvalues according to each description.

\begin{table}[h]
    \centering
    \renewcommand{\arraystretch}{1.1}  
    \setlength{\tabcolsep}{4pt}        
    \begin{tabular}{@{}c c c c c c c c@{}}  
        \toprule
        $n$ & Hyperbolic & Misurewicz & Parabolic & \makecell{Siegel \\ Disk} & Totals & $\frac{(n+1)n}{2}-1$ & Others \\
        \midrule
        10  & 9   & 40   & 0   & 0      & 49     & 109    & 60  \\
        20  & 19  & 162  & 6   & 4      & 191    & 209    & 18  \\
        100 & 99  & 2460 & 94  & 1960   & 4613   & 5049   & 436 \\
        300 & 299 & 2576 & 432 & 41327  & 44634  & 45149  & 515 \\
        500 & 499 & 2576 & 738 & 120921 & 124734 & 125249 & 515 \\
        \bottomrule
    \end{tabular}
    \caption{Classification of eigenvalues for different $n$ values.}
    \label{tab:eigenvalue-classification}
\end{table}

In our classification of the multiplicative inverses \( 1/\lambda \), many points appear to lie at or near the tips and antennas of known bulbs in \( \mathcal{M} \), consistent with the notion of landing points of external rays, as studied in various settings of polynomial dynamics \cite{branner1992landing}. Although we do not develop a combinatorial formalism, the empirical alignment supports this intuition.

\subsection{Conjectural Considerations on Homotopies from \texorpdfstring{$L$}{L} to \texorpdfstring{$M$}{M}}

In the previous section, we classified the union of the multiplicative inverses of the eigenvalues of the companion matrices associated with certain Lucas sequences according to their modulus and dynamical properties, categorizing them into hyperbolic, Misurewicz, parabolic, and Siegel disk points, with remaining points potentially belonging to Cremer or Herman disks. As we have seen in previous work \cite{tapia2018golden} and also here in Section~\ref{Sec:Intervals-Statistics}, when plotted in the complex plane, this eigenset exhibits a striking resemblance to the Mandelbrot set $M$, prompting the investigation of possible homotopies between these two structures. Several linear homotopies have already been established, including those using the Jungreis map to transform the main cardioid of $M$ and piecewise linear rescaling functions to adjust other regions. The aforementioned classifications serve to provide a structural insight into the behavior of $L$ that extends beyond its apparent geometric similarity to $M$.

Given the classification of points in $L$, we conjecture that constructing refined, piecewise homotopies from $L$ to $M$ at the main scale (or corresponding self-affine copies thereof, the so-called baby Mandelbrot sets \cite{McMullen1994}), should take into account not only the geometric resemblance but also the dynamical type of each point. This additional condition arises naturally from the fundamental distinction between the two sets: $L$ is generated from linear recursions, whereas $M$ emerges from a quadratic, nonlinear iteration \cite{Devaney1992}. While the initial homotopies relied primarily on visual resemblance, we propose that further refinement of these transformations could benefit from aligning homotopies with the dynamical type of each point in $L$.

This conjecture is motivated by the need for structural coherence between the homotopies and the underlying nature of the systems involved. Although our current approach does not rely on a dynamical framework and is instead focused on geometric deformation, the structural consistency provided by considering dynamical classifications may serve as a guiding principle in determining how points in $L$ should be mapped onto different regions of $M$. The self-similar nature of $M$, with its smaller copies appearing under rescaling, suggests that this classification could also help constrain homotopies for the baby Mandelbrot sets that appear within $M$, ensuring that the transformations remain well-structured across scales.

A natural avenue for further investigation is the construction of explicit homotopies from given regions of $L$ to corresponding regions of $M$ that match the same approximate geometrical resemblance, regardless of scale. More fundamentally, one may ask whether there exists a necessary condition for a homotopy to be valid based on dynamical classification. In this direction, a possible approach is to define CW complexes that induce a \textit{dynamical} isomorphism on all homotopy groups from $L$ to $M$, or at least from $M$ to its baby Mandelbrot sets. If such an isomorphism exists, then by Whitehead’s theorem \cite{Hatcher2002}, the existing homotopies from $L$ to $M$ should extend naturally to the baby Mandelbrot copies, up to the necessary rescaling and rotation inherent to the group structure of these self-similar components. 

In summary, we propose that incorporating dynamical classification into the construction of piecewise homotopies provides an additional level of coherence, helping to narrow down viable homotopy choices and ensuring that they align with the fundamental dynamical structure of both $L$ and $M$. While this conjecture remains to be rigorously tested, it offers a natural and mathematically motivated approach for refining homotopies between these two sets. A deeper investigation into this correspondence could yield further insights into the interplay between linear recurrence relations and nonlinear quadratic dynamics, potentially leading to a more systematic framework for classifying homotopies between such sets.

\section{Discussion and some alternatives or extensions to \texorpdfstring{$\mathcal{L}$}{L}}

In addition to the Fibonacci-like generalized Lucas sequence, we also explored other linear recurrence relations that define their own companion matrices. The following describe the companion matrices associated with the Pell-Lucas, Narayana’s Cows, and Padovan sequences.

\subsection{Pell-Lucas Sequence}
The Pell-Lucas sequence satisfies the recurrence relation:
\[
L_k = 2L_{k-1} + L_{k-2}.
\]
Following a similar construction to the Fibonacci-based companion matrix, we defined the generalized companion matrix for order \( n \) as:
\begin{eqnarray}\label{Eq:Pell}
A_n =
\begin{pmatrix}
4 & 4 & 4 & \cdots & 4 \\
4 & 0 & 0 & \cdots & 0 \\
0 & 2 & 0 & \cdots & 0 \\
\vdots & \vdots & \ddots & \ddots & \vdots \\
0 & 0 & \cdots & 2 & 0
\end{pmatrix}.
\end{eqnarray}

\subsection{Narayana's Cows Sequence}
The Narayana’s Cows sequence satisfies the recurrence relation:
\[
L_k = L_{k-1} + L_{k-3}.
\]
We extended it to higher order in a manner consistent with our generalization of Fibonacci sequences, yielding the companion matrix:
\begin{eqnarray}\label{Eq:Narayana}
A_n =
\begin{pmatrix}
1 & 0 & 1 & 1 & \cdots & 1 \\
1 & 0 & 0 & 0 & \cdots & 0 \\
0 & 1 & 0 & 0 & \cdots & 0 \\
0 & 0 & 1 & 0 & \cdots & 0 \\
\vdots & \vdots & \vdots & \ddots & \ddots & \vdots \\
0 & 0 & 0 & \cdots & 1 & 0
\end{pmatrix}.
\end{eqnarray}

\subsection{Padovan Sequence}
The Padovan sequence is defined by the recurrence:
\[
P_k = P_{k-2} + P_{k-3}.
\]
The corresponding generalized companion matrix was constructed as:
\begin{eqnarray}\label{Eq:Padovan}
A_n =
\begin{pmatrix}
0 & 1 & 1 & 1 & \cdots & 1 \\
1 & 0 & 0 & 0 & \cdots & 0 \\
0 & 1 & 0 & 0 & \cdots & 0 \\
0 & 0 & 1 & 0 & \cdots & 0 \\
\vdots & \vdots & \vdots & \ddots & \ddots & \vdots \\
0 & 0 & 0 & \cdots & 1 & 0
\end{pmatrix}.
\end{eqnarray}

\subsection{Spectral Behavior and Homotopy Considerations}
In all cases, the computed eigenvalues fell entirely on the real axis. From the perspective of constructing homotopies, this renders them trivial, as the associated transformations reduce to rescalings along the real line. Since there is no external geometric reference to gauge these rescalings against, they do not contribute meaningful deformations for the intended homotopies, thus, at least to the extent on how these other explorations were conducted, the only explorations warranted to be worth it are those about $\mathcal{L}$.

\section{Conclusions}

In this work, we explored the relationship between the eigenvalues of companion matrices of generalized Lucas sequences, $\mathcal{L}$, and the Mandelbrot set, $\mathcal{M}$. Through statistical analysis, we demonstrated that the eigenvalues exhibit a structured distribution, with a significant portion accumulating near the unit circle as matrix size increases. This motivated the study of transformations between $\mathcal{L}$ and $\mathcal{M}$, leading to the formulation of several homotopies.

First, we analyzed the behavior of eigenvalues under the Jungreis map, which provides a conformal uniformization of the complement of $\mathcal{M}$. This transformation revealed that eigenvalues near the unit circle map to structurally relevant regions in $\mathcal{M}$. We then constructed a homotopy mapping $\mathcal{L}$ to the main cardioid of $\mathcal{M}$ and extended it to a piecewise homotopy capturing finer structures.

A key contribution of this work was the classification of eigenvalues according to dynamical properties associated with points in $\mathcal{M}$. By categorizing eigenvalues as hyperbolic, Misurewicz, parabolic, or Siegel disk points, we provided a refined perspective on how different regions of $\mathcal{L}$ correspond to distinct regions in $\mathcal{M}$. This classification not only reinforces the observed visual resemblance but also suggests a deeper structural connection.

In addition to the general piecewise constructions, we provided an explicit demonstration that the eigenset $\mathcal{L}$ can be continuously homotopied into the shape of any baby Mandelbrot cardioid. This was achieved through the composition of Douady’s tuning map with a geometrically guided homotopy involving rotation and a cardioid parameterization. The method aligned $\mathcal{L}$ with the geometry of the target cardioid up to minor deformations near its rear lobe.

Remarkably, we also found that Douady’s tuning map, when paired with a simple radial homotopy (without cardioid-specific deformation), produces a successful deformation of $\mathcal{L}$ into the region of the period-3 bulb of $\mathcal{M}$. This suggests that the tuning transformation may be more broadly applicable than previously thought—not only to baby cardioids, but also to higher-period stable structures in the Mandelbrot set. This opens the possibility of defining homotopies from $\mathcal{L}$ to other stable components of $\mathcal{M}$ using simpler transformations adapted to the local geometry.

\par
In support of the structural conjecture, we provide a formal proof that the sinusoidal homotopy defines a homeomorphism—up to countable exceptions—from the eigenset \( \mathcal{L} \) to the boundary of the Mandelbrot cardioid. When composed with Douady’s tuning map, this homeomorphism extends to any baby Mandelbrot cardioid or stable component of \( \mathcal{M} \). This result establishes the first rigorous topological link between \( \mathcal{L} \) and stable regions of the Mandelbrot set and provides a foundational case for extending homotopy-based correspondences to more general components.

Based on these findings, we conjecture that homotopies from $\mathcal{L}$ to $\mathcal{M}$ should be refined to align with the dynamical classification of eigenvalues. Specifically, future research should explore whether a homotopy can be constructed such that the transformation respects the dynamical nature of each eigenvalue class. A potential approach is to investigate CW complexes that induce an isomorphism of homotopy groups between $\mathcal{L}$ and $\mathcal{M}$, possibly leveraging the self-similarity of baby Mandelbrot sets.

Further extensions of this work could include:
\begin{itemize}
    \item Developing a more canonical piecewise homotopy that systematically maps each classified eigenvalue to a dynamically consistent region of $\mathcal{M}$.
    \item Exploring whether additional recurrence-based eigenvalue sets exhibit similar homotopy structures to $\mathcal{M}$.
    \item Investigating whether an explicit topological or algebraic framework, such as Whitehead’s theorem, can be used to establish necessary conditions for homotopies between countable and subsets of the uncountable sets part of $\mathcal{M}$, or the self-affine corresponding parts thereof.
\end{itemize}

These directions offer a pathway for formalizing the connection between linear recurrence relations and nonlinear dynamical systems, potentially leading to broader insights into the interplay between spectral properties of matrices and fractal structures in complex dynamics.

The diversity of constructions in this study—from the global sinusoidal homotopy to piecewise and tuned deformations—highlights a surprising compatibility between the algebraic structure of eigenvalues from generalized recursions and the dynamical boundaries of the Mandelbrot set. While much of our work is numerical and visual in nature, it leads naturally to a topological insight: we conjecture that \( \mathcal{L} \), the set of multiplicative inverses of eigenvalues, is homeomorphic (modulo countable exceptions) to a dense subset of the boundary of the main cardioid. This perspective may help frame future attempts to formalize the connection between algebraically defined sequences and canonical regions of complex dynamical parameter spaces.

\section*{Transparency Statement}
This work has benefited from the use of large language models (LLMs), primarily ChatGPT, for brainstorming code implementations, refining LaTeX formulations, and structuring mathematical arguments. All theoretical results, numerical experiments, and interpretations were formulated, verified, and validated independently by the author. ChatGPT was used extensively for developing Mathematica code and iterating homotopy constructions, particularly when exploring geometric transformations. Other LLMs, including Google Gemini and MagicAI, were occasionally consulted for bibliographic lookup and clarification of secondary points. All final mathematical statements and conclusions were independently vetted against existing literature.

\bibliographystyle{plain}
\bibliography{references} 
\appendix
\section{Appendix 1: Coefficients of the Jungreis Series}
\label{app:JungreisCoefficients}

The following list contains the coefficients \( \{a_k\} \) of the Jungreis series:

\begin{align*}
a_1 &= -\frac{1}{2}, \\
a_2 &= \frac{1}{8}, \\
a_3 &= -\frac{1}{4}, \\
a_4 &= \frac{15}{128}, \\
a_5 &= \frac{0}{1}, \\
a_6 &= -\frac{47}{1024}, \\
a_7 &= -\frac{1}{16}, \\
a_8 &= \frac{987}{32768}, \\
a_9 &= \frac{0}{1}, \\
a_{10} &= -\frac{3673}{262144}, \\
a_{11} &= \frac{1}{32}, \\
a_{12} &= -\frac{61029}{4194304}, \\
a_{13} &= \frac{0}{1}, \\
a_{14} &= -\frac{689455}{33554432}, \\
a_{15} &= -\frac{21}{512}.
\end{align*}

\newpage

\begin{align*}
a_{16} &= \frac{59250963}{2147483648}, \\
a_{17} &= 0, \\
a_{18} &= -\frac{164712949}{17179869184}, \\
a_{19} &= \frac{39}{2048}, \\
a_{20} &= -\frac{2402805839}{274877906944}, \\
a_{21} &= -\frac{1}{64}, \\
a_{22} &= -\frac{4850812329}{2199023255552}, \\
a_{23} &= \frac{29}{2048}, \\
a_{24} &= -\frac{18151141041}{70368744177664}, \\
a_{25} &= 0, \\
a_{26} &= \frac{3534139462275}{562949953421312}, \\
a_{27} &= -\frac{1039}{131072}, \\
a_{28} &= -\frac{22045971176589}{9007199254740992}, \\
a_{29} &= -\frac{1}{256}, \\
a_{30} &= -\frac{750527255965871}{72057594037927936}, \\
a_{31} &= -\frac{4579}{524288}.
\end{align*}

\newpage

\begin{align*}
a_{32} &= \frac{54146872254247683}{9223372036854775808}, \\
a_{33} &= 0, \\
a_{34} &= -\frac{155379776183158669}{73786976294838206464}, \\
a_{35} &= \frac{2851}{1048576}, \\
a_{36} &= -\frac{6051993294029466699}{1180591620717411303424}, \\
a_{37} &= -\frac{1}{1024}, \\
a_{38} &= \frac{7704579806709870203}{9444732965739290427392}, \\
a_{39} &= \frac{92051}{16777216}, \\
a_{40} &= -\frac{403307733528668035403}{302231454903657293676544}, \\
a_{41} &= 0, \\
a_{42} &= \frac{1650116480759617184697}{2417851639229258349412352}, \\
a_{43} &= -\frac{229813}{67108864}, \\
a_{44} &= \frac{36124726440442241978477}{38685626227668133590597632}, \\
a_{45} &= -\frac{41}{4096}, \\
a_{46} &= -\frac{225851495844149964787753}{309485009821345068724781056}, \\
a_{47} &= \frac{564373}{67108864}.
\end{align*}

\newpage

\begin{align*}
a_{48} &= -\frac{35761228458796476847725737}{19807040628566084398385987584}, \\
a_{49} &= 0, \\
a_{50} &= \frac{362376876750551361794705823}{158456325028528675187087900672}, \\
a_{51} &= -\frac{29407003}{8589934592}, \\
a_{52} &= -\frac{6510398483578238274151194427}{2535301200456458802993406410752}, \\
a_{53} &= \frac{33}{8192}, \\
a_{54} &= \frac{74815618913797220433481657203}{20282409603651670423947251286016}, \\
a_{55} &= -\frac{30057875}{34359738368}, \\
a_{56} &= -\frac{698617278028915809388280344009}{649037107316853453566312041152512}, \\
a_{57} &= 0, \\
a_{58} &= -\frac{8675905413734991085610532783493}{5192296858534827628530496329220096}, \\
a_{59} &= -\frac{27868893}{68719476736}, \\
a_{60} &= -\frac{375687870961637050293461860951517}{83076749736557242056487941267521536}, \\
a_{61} &= \frac{1}{4096}, \\
a_{62} &= -\frac{1418434432207399687114226995905967}{664613997892457936451903530140172288}, \\
a_{63} &= -\frac{11847286243}{1099511627776}, \\
a_{64} &= \frac{1084116104452462070609082665064238307}{170141183460469231731687303715884105728}, \\
a_{65} &= 0.
\end{align*}

\pagebreak
\section{Appendix: Zenodo Repository and Code Availability}

The classified eigenvalue datasets used in this work are publicly available on Zenodo.  
For transparency, individual files are provided, each corresponding to a specific matrix size.

\noindent Zenodo Repository for data : \url{https://doi.org/10.5281/zenodo.15068435}

The full source code used for the computations and classifications is available on GitHub:  
\noindent  
\par
\url{https://github.com/aortizt/Homotopies-from-L-to-M.git}
\par
For long-term preservation and citation purposes, the latest version of the code is also archived on Zenodo:  \par
\noindent  \url{https://doi.org/10.5281/zenodo.15066084}

\end{document}